\documentclass{elsarticle}
\usepackage[utf8]{inputenc}
\usepackage[a4paper]{geometry}
\usepackage[colorlinks,linkcolor=black,urlcolor=black,citecolor=black]{hyperref}
\usepackage{amsmath, amsthm, amssymb, mathtools, tikz, enumitem, booktabs}
\usepackage[bf,singlelinecheck=false]{caption}
\usepackage[ruled]{algorithm}
\usepackage[noend]{algpseudocode}
\usepackage[english]{isodate}
\usepackage{algorithmicx, comment}
\usepackage{aliascnt}
\usepackage{multirow}
\usepackage{cleveref}
\usepackage{xspace}

\theoremstyle{definition}
\newtheorem{theorem}{Theorem}[section]
\newtheorem{definition}[theorem]{Definition}
\newtheorem{lemma}[theorem]{Lemma}

\newtheorem{remark}[theorem]{Remark}
\newtheorem{example}[theorem]{Example}

\newcommand{\actSet}{\textsc{actSet}\xspace}
\newcommand{\actEq}{\ensuremath{\preccurlyeq}}
\MakeRobust{\Call}

\def\<{\langle}
\def\>{\rangle}

\newcommand{\N}{\mathbb{N}}

\newcommand{\stackS}{S}
\newcommand{\stackT}{T}

\newcommand{\Refi}{\textsc{Ref}}
\newcommand{\Can}{\textsc{Canonical}}

\newcommand{\myfloor}[1]{\left \lfloor #1 \right \rfloor}

\newcommand{\approxFunc}{\textsc{Approx}}
\newcommand{\splitFunc}{\textsc{Split}}
\newcommand{\fixedFunc}{\textsc{Fixed}}
\newcommand{\partFunc}{\textsc{Part}}

\newcommand\Approx[2]{\ifthenelse{\equal{#2}{}}
                         {\textsc{Approx} (#1)}
                         {\Call{Approx}{#1,#2}}}
 \newcommand{\Split} {\textsc{Split}}
\newcommand{\Fixed}[1]{\Call{Fixed}{#1}}
\newcommand{\Part}[1]{\Call{Part}{#1}}

\newcommand{\Sym}[1]{\operatorname{Sym}\!\left(#1\right)}

\newcommand{\Auto}[1]{\operatorname{Aut}\!\left(#1\right)}
\newcommand{\Iso}[2]{\operatorname{Iso}\!\left(#1,#2\right)}
\newcommand{\Min}{\textrm{Min}}

\newcommand{\Stacks}[1]{\operatorname{\textsc{DigraphStacks}}\!\left(#1\right)}
\newcommand{\stacks}{\Stacks{\Omega}}

\newcommand{\LabelledDigraphs}[2]{\operatorname{\textsc{LabelledDigraphs}}\!\left(#1,#2\right)}
\newcommand{\EmptyStack}[1]{\operatorname{\textsc{EmptyStack}}\!\left(#1\right)}

\newcommand{\labelSet}{\mathfrak{L}}
\newcommand{\labelFunc}{\textsc{Label}}

\newcommand{\ignore}[1]{}

\begin{document}

\begin{frontmatter}

\title{Computing canonical images in permutation groups with Graph Backtracking}

\author{Christopher Jefferson}


\author
{Rebecca Waldecker\corref{mycorrespondingauthor}}
\cortext[rebecca.waldecker@mathematik.uni-halle.de]{Corresponding author}

\author{Wilf A. Wilson}

\begin{abstract}

    We describe a new algorithm for finding a canonical image of an object under the action of a finite permutation group. This algorithm builds on previous work using Graph Backtracking \cite{gb-published}, which extends Jeffrey Leon's Partition Backtrack framework \cite{leon1991,leon1997}. Our methods generalise both Nauty \cite{nauty}
    and Steve Linton's Minimal image algorithm~\cite{linton-smallest-image}.
\end{abstract}


\end{frontmatter}

\section{Introduction and background}\label{sec-intro}

Given a group \(G\) that acts on a finite set \(\Omega\), one of the simplest questions we can ask is whether two elements \(o_1,o_2 \in \Omega\) are in the same orbit of \(G\), and if so, to find an element \(g \in G\) such that \(o_1^g = o_2\).

This problem turns out to be surprisingly difficult. If we consider \(\Omega\) as the set of graphs on a vertex set \(V\), and \(G = \Sym V\), then the above question is equivalent to the graph isomorphism problem, which has been studied extensively both in theory \cite{graphisoquasi} and in practice. The best-known practical system for solving the graph isomorphism problem is Nauty \cite{nauty}.
The way in which Nauty solves this problem in practice, which we will also refer to in this article, is to use a \textbf{canonical labelling function}. Given a group \(G\) acting on a set \(\Omega\), a canonical labelling function maps each element \(\omega \in \Omega\) to a ``canonical'' member \(\Can(\omega)\) of its $G$-orbit, which we call its canonical image. This means that two elements of \(\Omega\) are in the same $G$-orbit if and only if they have the same canonical image, i.e. the same image under the canonical labelling function. 

Canonical labelling functions have several applications. As a brief example, suppose that we are given a set of graphs, all with vertex set $\Omega$, and we want to decide for any two given graphs whether or not they are in the same $\Sym\Omega$-orbit with respect to the naturally induced action. Given only two graphs, we could use a function which directly checks if the two graphs are in the same orbit, which may be faster than finding canonical images. However, for larger sets of graphs, if we can only check whether or not two graphs are in the same orbit, and if all the graphs are in different orbits, then we need to compare every graph to every other graph in order to prove that no two different graphs are in same orbit. If, instead, we can compute the canonical image of a graph under the action of $\Sym\Omega$, then we only need to look at each graph once, find its canonical image, and then we know that the graphs in the same orbit are exactly the ones with the same canonical image.

In previous work, we have discussed methods for calculating minimal and canonical images (see \cite{minimal-canonical}), building on the work of Linton (\cite{linton-smallest-image}).
Minimal images are a special case of canonical images, using a total ordering of the underlying set. This previous work was limited to the induced action of a permutation group on subsets of a set.
Another important special case is finding canonical images of graphs (Nauty, see \cite{nauty}). 
For completeness, we mention the related but very different problem of finding a canonical member of \textit{every} orbit under a given group action. This problem is treated in \cite{orbiter}.
There is also related work by Schweitzer and Wiebking on canonical images from a more general perspective (see \cite{WiebkingCanonical}), but we do not discuss it here. As far as we understand, their work considers the problem from a theoretical perspective and is not intended to lead to a practical algorithm, and in particular we cannot compare our methods to theirs in terms of performance. 

In this paper we describe a general and practical framework for finding canonical images. Our results build on Graph Backtracking (\cite{gb-published} and \cite{GBpkg}), which itself builds on Nauty and on the search method developed by Jeffrey Leon, known as Partition Backtracking. The use of the Graph Backtracking framework implies that we will re-use some concepts that we have introduced in previous articles, but with minor changes and with a stronger focus on practical implementation. While it is helpful to be familiar with \cite{gb-published}
when reading this article, it is not absolutely necessary, because we state the most important definitions and results and explain how we adapted them for this new work.

We begin with general comments and some examples for motivation, set-up and notation in Section~\ref{sec-prelim}.  
There we also adjust some earlier definitions
for the purpose of this article, and introduce some new definitions.
In Section 3 we explain the core of our technique, an how we build search trees, again with examples for illustration, and then in Section 4 we explain how we construct canonical images from such search trees.
Section 5 discusses some experiments, including comparisons with earlier work.
We close with a few final remarks and an appendix with some additional comments and results, including further details relevant for implementation.

\section{Preliminaries and first examples}\label{sec-prelim}

Throughout this paper \(\Omega\) is a finite set and \(G\) is a subgroup of \(\Sym\Omega\).
In our applications we often consider $G$ in its induced action on sets of lists, graphs with vertex set $\Omega$, subsets of $\Omega$, subgroups of $G$ (like point stabilisers) or other combinatorial structures.
We use the name \actSet for a set of combinatorial objects on which $G$ acts and we will find canonical images for members of \actSet in their orbit in $G$.
Throughout, we use standard notation for permutation groups and group actions without further explanation, in particular notation for orbits and for the stabilisers of points, lists or sets. See for example \cite{cameron}.

In order to efficiently determine whether two elements of \actSet are in the same orbit of \(G\), we apply the concept of a canonical image (see also Definition 2.1 from \cite{minimal-canonical}).

\begin{definition}
    A \textbf{canonical labelling function} \(C\) with respect to the action of \(G\) on a set \(\actSet\) is a function \(C:\actSet \mapsto \actSet\) such that, for all \(a \in \actSet\), it is true that:

    \begin{itemize}
        \item \(C(a) \in a^G\), and
        \item \(C(a^g) = C(a)\) for all \(g \in G\)
    \end{itemize}

    We call \(C(a)\) the \textbf{canonical image} of \(a\) (with respect to \(G\) and the action of \(G\) on \actSet).
    Whenever the group $G$ and its action on \actSet are clear from context, then we will not mention them.
    Further, we say that \(g \in G\) is a \textbf{canonising element} for \(a\) if and only if \(a^{g} = C(a)\).
\end{definition}

Inspired by work of Linton \cite{linton-smallest-image}, we introduced practical algorithms for minimal and canonical images in \cite{minimal-canonical}. In the present paper, we combine these ideas with the concept of Graph Backtrack search, and we give more explanations of implementation details.

One example of a canonical labelling function that is quite intuitive and easy to explain is given by the \textbf{minimal image}.
It is calculated with respect to a given total ordering $\actEq$ on \actSet, and it is defined as 
follows: 

For all $a\in \actSet$, 
\(\Min_{\actEq}(a) \coloneqq \min_{\actEq}\{a^g \mid g \in G\}\). 

Unfortunately, for a general canonical image, there is no concise and simple mathematical description of how they are created, because the resulting canonical images are affected by many minor parts of the algorithm.
But we can describe the idea in a nutshell:
We narrow down the possibilities for a canonical image, efficiently and respecting the action of $G$, and only at the very end we introduce an ordering and calculate a minimal image.
It is only in the last step where details like the exact orderings that are used on $\Omega$ and \actSet influence the outcome of the algorithm.
While even this might sound mildly concerning, or at least unsatisfying, this is what happens in practice all along, even for algorithms like Nauty. There, the exact implementation of equitable partitioning determines the canonical graph that is produced.   
In our algorithms, if the orderings chosen at the end are exactly the same, then the same canonical image will be calculated.
In contrast, a canonising element may not be unique.  More precisely, given a canonising element $g \in G$ for some object $a \in \actSet$, it follows that every element in the right coset $G_a \cdot g$ is a canonising element for $a$ as well.

We close the introduction to this section with two examples that have motivated our work.

\begin{example}\label{ex-motivation}

\mbox {}

(a) Nauty finds canonical images of undirected graphs. In our notation, with our given set $\Omega$, this means that $G=\Sym\Omega$ and we take \actSet to be the set of all undirected graphs with vertex set $\Omega$. Then there is a naturally induced action of $G$ on \actSet with well-defined orbits, and therefore we can apply our notion of a canonical labelling function.

(b) The canonical image algorithm in \cite{minimal-canonical}, which extends earlier work of Linton \cite{linton-smallest-image}, considers any subgroup \(G\le \Sym\Omega\), and \actSet is the set of all subsets of \(\Omega\). 
If $A \subseteq \Omega$, then $g \in G$ maps $A$ to \(A^g\coloneqq\{a^g \mid a \in A\}\).
With this induced action, we have well-defined $G$-orbits again and we can ask for canonical images of elements of \actSet.
\end{example}

While the two special cases from the previous example are important and occur frequently in practice, 
there are many other problems we might want to solve. 
As we will see later, many problems can be transformed and can then be solved as a classical ``graph problem'' or a classical ``canonical image problem''.
But there are limits to this approach and it does not always lead to practically useful algorithms, as can be seen in \cite{minimal-canonical}.

\subsection{Definitions and notation}

We keep this section as short as possible -- on the one hand we do not want to reproduce all the basic results and examples from previous papers, and on the other hand we do not expect all readers to have read \cite{minimal-canonical} and \cite{gb-published}, which is why we give a basic overview:

A \textbf{graph} with vertex set $\Omega$ is a pair $(\Omega, E)$, where $E$ is a set of $2$-subsets of
$\Omega$. A \textbf{directed graph} with vertex set $\Omega$, or \textbf{digraph} for
short, is a pair $(\Omega, A)$, where $A \subseteq \Omega \times \Omega$ is a
set of pairs of elements in $\Omega$ called \textbf{arcs}. The elements of
$\Omega$ are called \textbf{vertices} in the context of graphs and digraphs.  Our
definition allows a digraph to have \textbf{loops}, which are arcs of the form
$(\alpha, \alpha)$ for some vertex $\alpha \in \Omega$.

Our techniques for searching in $\Sym{\Omega}$ are built around digraphs in
which each vertex and arc is given a \textit{label} from a label set $\labelSet$.  We
define a \textbf{vertex- and arc-labelled digraph}, or \textbf{labelled digraph} for
short, to be a triple $(\Omega, A, \labelFunc)$, where $(\Omega, A)$ is a
digraph and $\labelFunc$ is a function from $\Omega \cup A$ to $\labelSet$.
We often refer to labels for vertices as colours and say that vertices are coloured or individualized. 

Let $\LabelledDigraphs{\Omega}{\labelSet}$ denote the class of all
labelled digraphs on $\Omega$ with labels in $\labelSet$.
A \textbf{labelled digraph stack} on $\Omega$ is a finite (possibly empty) list of
labelled digraphs on $\Omega$.  We denote the collection of all labelled digraph
stacks on $\Omega$ by $\Stacks{\Omega}$.  The \textbf{length} of a labelled
digraph stack $\stackS$, written $|\stackS|$, is the number of entries that it
contains.  A labelled digraph stack of length $0$ is called \textbf{empty}, and we
denote the empty labelled digraph stack on $\Omega$ by $\EmptyStack{\Omega}$.
We use a notation typical for lists, whereby if $i \in \{1, \ldots,
|\stackS|\}$, then $\stackS[i]$ denotes the $i^{\text{th}}$ labelled digraph in
the stack $\stackS$.
If  $\stackS, \stackT \in \Stacks{\Omega}$, then we write 
  $\stackS \cong \stackT$ if and only if
    $\stackS$ and $\stackT$ are \textbf{isomorphic}, which means that they have the same length and there exists a permutation $g \in \Sym{\Omega}$ such that for all
    $i \in \{1, \ldots,
|\stackS|\}$, $\stackS[i]^g = \stackT[i]$. The set of isomorphisms from $\stackS$ to $\stackT$ is denoted by $\Iso{\stackS}{\stackT}$, and along the same lines, the automorphism group of the stack $\stackS$ is denoted by $\Auto{\stackS}$.

We allow any labelled digraph stack on $\Omega$ to be appended
onto the end of another.  If $\stackS, \stackT \in
\Stacks{\Omega}$ have lengths $k$ and $l$, respectively, then we define $\stackS
\Vert \stackT$ to be the labelled digraph stack
$\left[
    \stackS[1], \ldots, \stackS[k],
    \stackT[1], \ldots, \stackT[l]
    \right]$
of length $k + l$ formed by appending $\stackT$ to $\stackS$.

Now we discuss some concepts that have appeared in earlier work, but that we will refine and specify for
the algorithms relevant to this article.
The main reason we have to adjust the existing definitions is that the proofs in \cite{gb-published} and \cite{gb-extended} are only concerned with the behaviour of a single search -- in standard Graph Backtracking a single problem is solved, and a group, coset, or the empty set is found.  
When finding canonical images, we must prove that searches for the canonical image of two elements of \actSet in the same orbit of $G$ produce the same canonical image and therefore, in several places, we need to add extra conditions about the 
compatibility with the action of $G$ or $\Sym\Omega$.

The following definitions are in the spirit of \cite{gb-published}, but they are often adapted for the purpose of this article. There will be several examples in \Cref{sec-approx} and later.

\begin{definition}\label{defn-comp-refiner}
If \(a \in \actSet\), then a \textbf{refiner for \(a\)} is defined to be a function 
$\Refi$
 from $\Stacks{\Omega}$ to itself such that, whenever
  $\stackS, \stackT \in \Stacks{\Omega}$ and
  $\stackS \cong \stackT$, it is true that
  \[
    \Sym\Omega_a \cap \Iso{\stackS}{\stackT}
    \subseteq
    \Sym\Omega_a \cap \Iso{\Refi(\stackS)}{\Refi(\stackT)}.
  \]

We say that a refiner $\Refi$ for an object $a \in \actSet$ is \textbf{compatible with} $\Sym\Omega$ if and only if for all  
$\stackS \in \Stacks{\Omega}$ and all \(g \in \Sym{\Omega}\), we have that
\(\Refi(\stackS^g) = (\Refi(\stackS))^g\).  
\end{definition}

The main difference between Definition \ref{defn-comp-refiner} and the refiners from \cite{gb-published} is that we now define refiners for elements \(a\) of \(\actSet\), and then use the group $\Sym\Omega_a$, rather than restricting to subgroups of $\Sym\Omega$ directly. The additional compatibility condition is necessary here -- it captures the fact that refiners for two objects that are in the same orbit must have refiners which operate in ``the same way'' -- this will ensure we find the same canonical image for the two objects. Compatibility will therefore be used in the proofs of the correctness of our algorithms in Sections 3 and 4.

Often, compatibility is automatically fulfilled, as is the case for most of our previously published refiners. But since we did not introduce this notion in earlier work (because it was not needed), we did not check our refiners in \cite{gb-published} for compatibility and therefore we discuss it here. It will turn out that several related concepts must be slightly adapted as well.

When performing our backtracking search, we need to find a way of splitting the current search node into multiple child nodes. In Graph Backtracking, the function used to split has relatively few restrictions, apart from ensuring that the entire space of possible solutions is partitioned between the new child nodes and each child represents a strictly smaller part of the search than its parent.

To ensure that we find the same canonical image when looking at different objects with the same canonical image, we need to more tightly constrain how we split search nodes, taking the action of the group $\Sym\Omega$
on \actSet into account. Once more we adapt definitions from \cite{gb-published}:

\begin{definition}\label{def-approx-iso}
  A \textbf{compatible isomorphism approximator} is a function
  $\approxFunc$ which maps a pair of labelled digraph stacks on $\Omega$ to a subgroup of $\Sym{\Omega}$, a coset of a subgroup of $\Sym{\Omega}$, or the empty set, in such a way that the following conditions are satisfied
  for all $\stackS, \stackT \in
  \Stacks{\Omega}$ and $g \in \Sym{\Omega}$:
 
  \begin{enumerate}[label=\textrm{(\roman*)}]
    \item\label{item-approx-true-overestimate}
      $\Iso{\stackS}{\stackT} \subseteq \approxFunc(\stackS,\stackT)$.

    \item\label{item-approx-different-lengths}
      If $|\stackS| \neq |\stackT|$, then
      $\approxFunc(\stackS,\stackT) = \varnothing$.

    \item\label{item-approx-right-coset-of-aut}
      If $\approxFunc(\stackS,\stackT) \neq \varnothing$, then there is some $h \in \Sym{\Omega}$ such that
      $\approxFunc(\stackS,\stackT) = \approxFunc(\stackS,\stackS) \cdot h$.
      
     \item $\approxFunc(\stackS^g,\stackT^g)=\approxFunc(\stackS,\stackT)^g$.
  \end{enumerate}
  
Given $\stackS \in
  \Stacks{\Omega}$, we usually abbreviate $\approxFunc(\stackS,\stackS)$ as $\approxFunc(\stackS)$. 
\end{definition}

The isomorphism approximators defined in Section 5 of \cite{gb-published} are all compatible already. The first, given in Section 5.1, calculates $\Iso{\stackS}{\stackT}$ exactly, hence it is naturally compatible: For all graph stacks \(\stackS\) and \(\stackT\) and all permutations $g \in \Sym\Omega$, it is true that $\Iso{\stackS^g}{\stackT^g} = \Iso{\stackS}{\stackT}^g$.

Another way of viewing Condition (iv) of \Cref{def-approx-iso} is that the isomorphism approximator commutes with re-naming the elements of \(\Omega\). Using this observation, we can see that the other approximators from Section 5 of \cite{gb-published} are also compatible, because none of the approximations ever refer to an ordering, or particular elements, of \(\Omega\). We next similarly extend the definition of a splitter from Section 6 of \cite{gb-published}, again adding an extra requirement to ensure compatibility.

\begin{definition}\label{def-invariant-splitter}
Suppose that $\approxFunc$ is a compatible isomorphism approximator and 
let $\mathcal{T}$ denote the set of labelled digraph stacks $T$ on $\Omega$ where $|\approxFunc(T)| > 1$.
Let $\Split$ be a function from $\mathcal{T}$ to the set of finite ordered lists of labelled digraph stacks on $\Omega$.
  
  We say that $\Split$ is a \textbf{compatible splitter} (with respect to $\approxFunc$) 
  if and only if, for all $T \in \mathcal{T}$ and $n \in \N$, the ordered list $\Split(T)=[T_1,...,T_n]$ satisfies the following:

  \begin{enumerate}[label=\textrm{(\roman*)}]
    \item
      $\Iso{\stackT}{\stackT}
        =
        \Iso{\stackT \mathop{\Vert} \stackT_{1}}{\stackT \mathop{\Vert} \stackT_{1}}
        \mathop{\cup}
        \Iso{\stackT \mathop{\Vert} \stackT_{1}}{\stackT \mathop{\Vert} \stackT_{2}}
        \mathop{\cup}
        \cdots
        \mathop{\cup}
        \Iso{\stackT \mathop{\Vert} \stackT_{1}}{\stackT \mathop{\Vert} \stackT_{n}}$.

      \item
      $|\Approx{\stackT \mathop{\Vert} \stackT_{1}}{\stackT \mathop{\Vert} \stackT_{i}}|
       <
       |\approxFunc(T)|$ for all $i \in \{1, \ldots, n\}$.

    \item
      For all \(g \in \Sym\Omega\), the sets $\{T_1,...,T_n\}^g$ and 
$\{T_1^g,...,T_n^g\}$ coincide. This means that the lists \(\Split(\stackT^g)\) and \((\Split(\stackT))^g\) are equal up to ordering. 
       
  \end{enumerate}

 
\end{definition}

How can we create compatible splitters? The ``fixed point splitter'' of~\cite[Definition~6.4]{gb-published}, does not satisfy our requirements, because it takes the smallest member of \(\Omega\) which is in a shortest orbit (ignoring fixed points). This will usually not be compatible with the application of permutations. In the next subsection, we build some necessary machinery to create compatible splitters.

\subsection{Approximators for canonisation}\label{sec-approx}

To create compatible splitters, we first introduce a way of ordering (our current approximation of) the orbits of the group. 

\begin{definition}\label{def-ooa}

An \textbf{ordered orbit estimator} for labelled digraph stacks is a function \partFunc{} that maps 
a labelled digraph stack on \(\Omega\) to 
an ordered partition of \(\Omega\) such that for each \(S \in \Stacks{\Omega}\) the following holds:

\begin{itemize}
    \item If \( i,j \in \Omega\) are in the same orbit of \(\Auto{S}\), then they are in the same cell of \Part{S}, and
    \item \(\Part{S}^g = \Part{S^g}\) for all \(g \in \Sym{\Omega}\).
\end{itemize}

From there we can define an 
\textbf{ordered orbit approximator} \approxFunc$_O$ as follows:
Given two labelled digraph stacks \(S,T \in \Stacks{\Omega}\), we define 
\approxFunc$_O(\stackS,\stackT)$ as the set of permutations in $\Sym\Omega$ that map \Part{S} to \Part{T}.
As before, if $\stackS \in
  \Stacks{\Omega}$, then we abbreviate \approxFunc$_O(\stackS,\stackS)$ as \approxFunc$_O(\stackS)$. 
\end{definition}

The first property in our definition of an ordered orbit estimator can be expressed as "\(\Part{S}\) is (possibly) coarser than the orbit partition of \(\Auto{S}\)". The reason why we allow a coarser partition is that calculating the exact orbits of \(\Auto{S}\) would require solving a graph isomorphism problem.
For the corresponding approximator, we note that just as for isomorphism approximators, the images under this map are subgroups of $\Sym\Omega$, cosets of a subgroup or the empty set.
Next we give examples and we show how an ordered orbit estimator may merge together different \(\Auto{S}\)-orbits:

\begin{example}\label{ex-ord-orbit-approx}
Suppose that $\Omega\coloneqq\{1,2,3,4,5,6,7,8\}$ and that $\Gamma_1$ is the graph with vertex set $\Omega$ and the following edges (viewed as arcs both ways):
$\{1,2\},\{1,3\},\{2,3\}, \{5,6\},\{6,7\}$.
Let $S$ be the digraph stack with just the entry
$\Gamma_1$.

Also, let $\Gamma_2$ be the graph with vertex set $\Omega$ and edges
$\{1,6\}, \{1,7\}, \{2,6\}, \{2,7\}, \{3,8\}, \{4,8\}$, and $\{5,8\}$, again viewed as arcs both ways.
Let $T$ be the digraph stack with just the entry
$\Gamma_2$.

\smallskip
(a) Let $F$ denote the function that, for each graph in the stack, creates an ordered partition where the cells correspond to the number of neighbours that each vertex has. In practice, we need to decide in which order to organise the cells and apply this ordering
consistently.
For example, $S$ could get mapped to
$F(S)\coloneqq[4,8|5,7|1,2,3,6]$
and $T$ could get mapped to 
$F(T)\coloneqq[3,4,5|1,2,6,7|8]$.
Just to clarify our notation for ordered partitions: The vertical lines seperate the cells, and the cells themselves are considered to be sets (hence unordered), while the ordering of the cells within the partition matters. 

In order to check, at least for these two stacks, 
whether or not $F$ is an ordered orbit estimator in the sense of Definition \ref{def-ooa}, we need to
consider the orbits of $\Auto S$, namely
$\{1,2,3\}, \{4,8\}, \{5,7\}$, and $\{6\}$,
and the orbits of $\Auto T$, which are 
$\{1,2,6,7\}, \{3,4,5\}$, and $\{8\}$.

A simple check shows that all elements that are in the same orbit under the automorphism group also are in the same cell. 
Also, the number of neighbours of vertices is invariant under the action of $\Sym\Omega$, so that both conditions for an ordered orbit estimator are satisfied. The image \approxFunc$_O(\stackS,\stackT)$ under the corresponding approximator is the empty set, because there is no element in $\Sym\Omega$ that maps $F(S)$ to $F(T)$.

\smallskip
(b) Another possibility would be the function $F$ that, for each graph in the stack, creates an ordered partition where the cells correspond to connected components.
Again, in practice, we need to decide how to order the cells.
For example $F(S)=[4|8|1,2,3|5,6,7]$ and $F(T)=[1,2,6,7|3,4,5,8]$.
We do our simple check, and this time we do not obtain an ordered orbit estimator.
The isolated vertices 4 and 8 in the graph $\Gamma_1$ are in the same $\Auto S$-orbit, but in different cells.
For $\Gamma_2$, we see that the vertex 8 is not in the same $\Auto T$-orbit as 3, 4 and 5, but it is in the same cell
of $F(T)$. Therefore, taking the connected components does not create a valid ordered orbit estimator.
However, taking connected components commutes with the action of $\Sym\Omega$, it is therefore not an unnatural way to produce a partition. 

Rephrasing the first condition of an ordered orbit estimator gives another way to look at this:
If  $\stackS \in
  \Stacks{\Omega}$, then
$\Auto S$ is a (possibly proper) subgroup of
the stabiliser of the ordered partition $\Part S$ in $\Sym\Omega$.
In that light,  connected components produce too fine a partition, because members of the same orbit of the automorphism can be in different connected components. While the automorphism group can permute connected components, it cannot permute its own orbits non-trivially. 
\end{example}

\smallskip
Using ordered orbit estimators and approximators, we can create a compatible splitter, as demonstrated in the next example:

\begin{example}[Compatible splitter]\label{ex-equitable-splitter}

\mbox{}

Suppose that \(E\) is an ordered orbit estimator on
\(\stacks\) and that 
\approxFunc$_E$ denotes the corresponding approximator as in Definition \ref{ex-ord-orbit-approx}.
Now we define a compatible splitter \(\textsc{Split}_E\) as follows: For any 
$\stackS \in \stacks$, 
calculate $E(S)$.
If all cells are singletons, then we do not consider $S$ because there is a unique permutation in $\Sym\Omega$ that maps $E(S)$ to itself, and hence
$|E(S)|=1$.

Otherwise we determine the size $k_S$ of the shortest 
non-singleton cell and we take the first 
cell in $E(S)$ of this minimal size $k_S$, and then generate a list $[S_1,...,S_{k_S}]$ of graph stacks, each consisting of just one graph.
For
every element in this shortest non-singleton cell, this unique graph in the stack is such that the vertex for this element is individualised.  This means the graph has no edges, all vertices except the individualised vertex are labelled $0$, and the individualised vertex is labelled $1$. (The exact labels are unimportant, but the same labels must always be used). The automorphism group of this graph fixes the individualised vertex and allows all other vertices to freely permute.
Finally, we set \(\textsc{Split}_E(S):=[S_1,...,S_{k_S}]\) as described above.

Going back to Example \ref{ex-ord-orbit-approx}, we make this explicit and look at 
the ordered orbit estimator $F$ that we introduced in (a). The images were
$F(S)=[4,8|5,7|1,2,3,6]$
and
$F(T)=[3,4,5|1,2,6,7|8]$.

For each number $i \in \{1,...,8\}$, let $\Delta_i$ denote the stack that consists solely of the empty graph 
on $\{1,...,8\}$ where the vertex $i$ is individualised as described above.
Then the ordering of the cells that was used for $F$ gives that 
\(\textsc{Split}_F(S)=[\Delta_4,\Delta_8]\), because $2$ is the size of the shortest non-singleton cell and the first such cell is $\{4,8\}$,
and
\(\textsc{Split}_F(T)=[\Delta_3,\Delta_4, \Delta_5]\), because here we have a unique shortest non-singleton cell, namely $\{3,4,5\}$.

\end{example}

In some situations we will need a fixed point approximator, which was discussed in our previous work and which gives a subset of the fixed points set of the graph stack. This can be easily created from an ordered partition approximator:

\begin{definition}
We define the function \textsc{Fix} from the set of all ordered partitions 
on \(\Omega\) to the set of all ordered lists with entries from \(\Omega\) as follows:

For each ordered partition \(P\) of $\Omega$, 
$\textsc{Fix}(P)$ is defined to be the list of elements of $P$ that are in singleton cells, in the order in which they appear in $P$.
\end{definition}

As an example for this last definition, we look at $\Omega\coloneqq\{1,...,8\}$ and we consider $P\coloneqq[4|8|1,2,3|6|5,7]$. 
Then $\textsc{Fix}(P)=[4,8,6]$.

In line with the compatibility conditions that we have introduced so far, we introduce a compatible partial order, which we will use for ordering the nodes of our search trees.

\begin{definition}
A \textbf{compatible partial order on the set of digraph stacks on \(\Omega\)} is a partial order $\preceq$ on  $\stacks$ such that,
for all \(S, T \in \stacks\) and all \(g_1,g_2 \in \Sym{\Omega}\): 
\[S \preceq T \text{~if and only if~} S^{g_1} \preceq T^{g_2}.\]
\end{definition}

For the next definition we recall that the images of an ordered orbit estimator are ordered partitions of $\Omega$, and we give an example of a compatible partial order. It will be discussed in \Cref{lem-comp-part-order}.

\begin{definition}\label{def-comp-part-order}
Suppose that \(\partFunc\) is an ordered orbit estimator on \(\stacks\) and let $S,T \in \stacks$.

We define \(S \preceq T\) if and only if 
one of the following holds:
The list \([~|O| \mid O \in \Part{S}~]\) is shorter than the list \([~|O| \mid O \in \Part{T}~]\), or  the lists have the same length and \([~|O| \mid O \in \Part{S}~]\) comes first under lexicographic ordering.
\end{definition}

\begin{example}\label{ex-comp-part-order}
We build on Example \ref{ex-ord-orbit-approx}, and again we let 
$\Omega\coloneqq\{1,2,3,4,5,6,7,8\}$ and we use the graphs $\Gamma_1$ and $\Gamma_2$ as well as the stacks $\stackS$ and $\stackT$ from that example, along with the function 
$F$ that, for each graph in the stack, creates an ordered partition where the cells correspond to the number of neighbours that each vertex has. We recall that
$F(S)=[4|8|5,7|1,2,3,6]$ and
$F(T)=[3,4,5|1,2,6,7|8]$ once we decided that we begin with the lowest number of neighbours.

In order to compare $\stackS$ and $\stackT$ with $\preceq$ from the previous definition, we write down the lists of cell sizes.
For $F(S)$, we obtain $[1,1,2,4]$, and for
$F(T)$, we obtain $[3,4,1]$, which means that
$T \preceq S$.
Now we check for compatibility, hence we let $g_1,g_2\in \Sym\Omega$.
Given the action of $\Sym\Omega$ on $\LabelledDigraphs{\Omega}{\labelSet}$, it follows that $F(\stackS^{g_1})=F(\stackS)$ if we 
order the cells in the same way, starting with the lowest number of neighbours.
Also $F(\stackT^{g_2})=F(\stackT)$, for the same reason, and consequently
$T^{g_2} \preceq S^{g_1}$.
Therefore we have an example of a compatible partial order, and this is in fact the ordering we will use later on.
\end{example}

The following lemma shows that the construction from the previous example works more generally:

\begin{lemma}\label{lem-comp-part-order}
Let \(\partFunc\) be an ordered orbit estimator on \(\stacks\) and let \( \preceq \) be as defined in \Cref{def-comp-part-order}.
Then this order is a compatible partial order.
\end{lemma}

\begin{proof}
We have seen the main argument in Example
\ref{ex-comp-part-order}.
Let $S,T \in \stacks$ be such that $S \preceq T$ and let $g_1,g_2\in \Sym\Omega$.

The compatibility of \Part{} implies that \(\Part{S}^{g_1} = \Part{S^{g_1}}\), and therefore 
\([~|O| \mid O \in \Part{S^{g_1}}~]\) = \([~|O| \mid O \in \Part{S}^{g_1}~]\) = \([~|O| \mid O \in \Part{S}~]\), because the action of permutations on ordered partitions respects the ordering and size of the cells.
For the same reason 
\([~|O| \mid O \in \Part{T^{g_2}}~]\) = \([~|O| \mid O \in \Part{T}~]\), which implies that
$S^{g_1} \preceq T^{g_2}$.

Since permutations are bijective, we can argue in the same way for the reverse implication.
\end{proof}

\section{Searching with graph stacks}

Now that we have introduced compatible refiners and compatible splitters, we will show how we use these together to search for candidates for canonical images. This search is split into two parts. The present section will build most of the machinery, using a backtrack search with graph stacks and then return a list of candidates. The next section will show how to systematically pick a canonical image from this list of candidates.
While this algorithm uses much of the same machinery as Graph Backtracking, it behaves quite differently. Most notably, in Graph Backtracking there is always a pair of graph stacks, called the left and right hand stacks, and the search stops if there are no isomorphisms from the left hand stack to the right. When looking for canonical images, we will only consider a \emph{single} graph stack. This stack corresponds to the \emph{right hand side stack in Graph Backtracking}.

\subsection{Building search trees}

Suppose that, for some finite set \(\Omega\), some group \(G \leq \Sym{\Omega}\) and some set of objects \actSet, we wish to find canonical images for members of \actSet with respect to the action of $G$. Our algorithm requires the following: a set of compatible refiners for the elements of \actSet that we wish to find canonical images for (only one of these will be used during any invocation of the algorithm), a compatible splitter \(\splitFunc\), an ordered orbit estimator  \partFunc{} and corresponding approximator \approxFunc,
a refiner \(R_G\) for \(G\). For $R_G$  we use the original definition of a refiner from \cite{gb-published}, 
because $R_G$ does not need to be compatible. There are several different choices for $R_G$, for example we can use the orbits of $G$. 
We will discuss more options in the appendix of this paper.

All the concepts mentioned above will be used throughout this section without reintroduction. Note that changing any of the refiners, the splitter, or the ordered orbit estimator and approximator can change the resulting canonical image.


We begin by giving a slightly simplified explanation of our method in \Cref{alg-basic-search}. Note that this initial presentation shows an algorithm which correctly returns a set of permutations, namely a set of candidates for 
permutations that map  \(a \in \actSet\) to a canonical image.
However, this set is produced in a very inefficient way. We will then show how the search can be improved. This second algorithm uses the function \textsc{MinPerm} which, given a group \(G\) and a list of integers \(L\), returns a permutation \(g \in G\) which maps \(L\) point-wise to its lexicographic minimal image. An implementation of this function is given in the Appendix.

\begin{algorithm}[!ht]
  \caption{A recursive algorithm using labelled digraph stacks to search
    in $\Sym{\Omega}$.}\label{alg-basic-search}

  \begin{algorithmic}[1]

    \item[\textbf{Input:}] \(a \in \actSet\) and a compatible refiner \(R_a\) for \(a\), taken from a set of compatible refiners for \(\actSet\).

    \item[\textbf{Output:}]
    A set of permutations \(g \in G\) containing the permutations which map \(a\) to its canonical representative.

    \vspace{2mm}
    \Procedure{Search}{$\stackS$}
    \Comment{The main recursive search procedure}\label{line-search-proc-start}

    \State{\(\stackS' \gets \Call{Refine}{\stackS}\)}
    \Comment{Refine the stack}\label{line-refine}

    \If{$|\approxFunc(\stackS')|=1$}
    \label{line-if-case-one-candidate}
      \State{\Return{\textsc{MinPerm}$(G,\textsc{Fix}(\partFunc(\stackS')))$}}  \label{line-case-one-candidate}
    \Else{}      \Comment{Multiple potential solutions}\label{line-multiple-candidates}

      \State{\([S_1, \dots, S_t] \gets \Call{Split}{\stackS'}\)}
      \Comment{Split the search space}\label{line-split}

      \State{\Return{\(
        \displaystyle\bigcup_{i \in \{1,\ldots,t\}}
        \Call{Search}{\stackS' \mathop{\Vert} \stackS_{i}}\)}}
      \Comment{Search recursively}\label{line-recurse}
    \EndIf{}

    \EndProcedure{}

    \vspace{1mm}
    \Procedure{Refine}{$\stackS$}
    \Comment{Attempt to prune the search space}\label{line-refine-proc-start}

    \While{\textsc{True}}
    \Comment{Loop until fixed point}\label{line-while-begin}\label{line-refine-empty}

    \State{\(\stackS_C \gets \stackS\)}
    \Comment{Save the stack to check if progress is made}\label{line-refine-store}

    \For{\(R \in [R_a, R_G]\)}\label{line-refine-loop}
      \State{$\stackS \gets \stackS \mathop{\Vert} R(\stackS)$}
      \Comment{Apply each refiner in turn}\label{line-apply-refiner}
    \EndFor{}

    \If{\(|\approxFunc(\stackS)| \not<
      |\approxFunc(\stackS_C)|\)}\label{line-refine-not-smaller}

      \State{\Return{\(\stackS\)}}
      \Comment{Stop: the refinement is no longer making progress}\label{line-refine-loop-end}

    \EndIf{}

    \EndWhile{}

    \EndProcedure{}

    \vspace{5mm}
    \State{\Return{\Call{Search}{$\EmptyStack{\Omega}$}}}\label{line-start-search} \Comment{Start the search}

  \end{algorithmic}
\end{algorithm}

We can view an execution of \Cref{alg-basic-search} as a search tree (\Cref{def-tree}).

\begin{definition}\label{def-tree}
A finite directed graph is a \textbf{search tree} if and only if it has the following properties:

\begin{itemize}
\item
It is connected,

\item
it has a unique initial node (the \textbf{root node}) with no incoming arcs and only outgoing arcs, 

\item 
every node except for the root node has exactly one incoming arc, and

\item
the graph contains no cycles (neither directed nor undirected).
\end{itemize}
This means that for every node, there is a unique directed path in the graph from the root node to it. The \textbf{children} of a node \(v\) are all nodes which are connected to an outgoing arc of \(v\). The nodes that only have an incoming arc and no outgoing arcs will be called \textbf{leaves}. 

We allow a search tree to be vertex-labelled and we place no restrictions on labels used on vertices.
\end{definition}

We now describe \Cref{alg-basic-search} in terms of a search tree by creating a node for each call of \textsc{Search(S)}, labelling this node with the \emph{graph stack} \(S\). We define a directed edge from any node to the nodes which are created by the recursive calls to \textsc{Search} \Cref{line-recurse} from that node.  Leaf nodes (which originate from calls to \textsc{Search(S)} which return on \Cref{line-case-one-candidate}) receive an additional label, namely a \textbf{result} permutation, which is the permutation returned on \Cref{line-case-one-candidate}.


In the next section we describe how to find a canonical image from the search tree. We will not reproduce all results and definitions from our previous work, for more technical details see \cite{gb-published} and \cite{perfect-refiners}. We now give a full example of \Cref{alg-basic-search}. This example does not demonstrate all possible cases. For a more complete discussion of applying refiners during search, and the issues involved, see \cite{gb-published}.

\begin{example}\label{example:canonical-image-graph}
\ \\
In this example,  $\Omega\coloneqq\{1,...,8\}$, \(G \coloneqq \langle (1,2,3,4,5,6,8), (1,3,2,6,4,5),(1,6)(2,3)(4,5)(7,8)\rangle\le \Sym\Omega\)
and our object \(a\) from \actSet is a graph $\Gamma$ with vertex set $\Omega$ and list of edges $\{1,6\}, \{1,7\}, \{2,6\}, \{2,7\}, \{3,8\}, \{4,8\}, \{5,8\}$. For our approximator we use \textsc{Approx$_S$} from Definition 5.9 of \cite{gb-published}, and we explain explicitly what it does in this specific example.
The refiners that we use are \(R_a\) and \(R_G\). For
\(R_a\), we convert $\Gamma$ into a digraph by replacing each edge with arcs in both directions, 
and then we add this graph to the stack (if it is not already in the stack). The refiner \(R_G\) adds a graph that connects vertices, again with arcs in both directions, if and only if they are in the same $G$-orbit. 

Now we can build a search tree for $\Gamma$ 
as described in \Cref{alg-basic-search}.
At the first call of \textsc{Search} (which is also the root node), \(S\) is the empty stack. We begin by calling \textsc{Refine}. We call \(R_a\) and \(R_G\), and as explained above this means that we add a directed version of $\Gamma$ and then a graph whose connected components are exactly the $G$-orbits.

We then calculate the image under the ordered orbit approximator \textsc{Approx}, as we have seen in \Cref{ex-ord-orbit-approx}~(a).
There are three vertices in the graph \(a\) with exactly one neighbour (namely 3,4 and 5), four vertices with two neighbours (namely 1,2,6,7) and there is a unique vertex with three neighbours (namely 8). For our bookkeeping in this example we write the image down as an ordered partition of $\Omega$: $P_1 = [3,4,5|1,2,6,7|8]$. Note that the refiner for \(R_a\) does not add anything anymore, because the graph stack already contains a directed version of the original graph \(a\)), so we will not discuss it further in this example.

We call \(R_G\) again, using the fact that $8$ is the only singleton cell, so we can consider the orbits of $G_8$. This does not produce any new information, because $G$ is $2$-transitive. Thus we do not change the result of \textsc{Approx}, and therefore \textsc{Refine} returns.
Since we do not (yet) have a partition where all cells are singletons (see \Cref{line-if-case-one-candidate} of the algorithm), we move on to calling \textsc{Split}. The splitter uses the partition \(P_1\) and branches on one of its cells. In theory, the ordering of the values returned by \textsc{Split} does not matter, but it can make a significant difference in terms of practical performance.

We choose a shortest non-singleton cell in the partition, which at this stage means that we split at the cell $[3,4,5]$. This results in three graphs, each with one vertex coloured
(individualising the corresponding point). Each of these graphs gives a child node, and a recursive call to \textsc{Search}.

Down each branch, we go through \textsc{Refine} again. Since $G$ is $3$-transitive, the refiners for these point stabilisers do not give graphs which change the result of \textsc{Approx}. 
Continuing down the branch for the vertex $3$, we call \textsc{Split} again and colour $4$ or $5$, which always individualises both vertices (just ordered differently).

After this splitting, we have now fixed the list of vertices $[8,3,4,5]$ (or $[8,3,5,4]$ in the other branch). $G_{[8,3,4,5]}$ is the trivial group, and therefore the refiner $R_G$ can now give every vertex a different colour. We represent this by ordering the vertices. However, the refiner must still be compatible, meaning this order must be compatible between different branches.

In our implementation, we choose an ordering for the orbits by using the minimal image of the list $[8,3,5,4]$ of fixed points under the action of $G$.
In this case, we find a permutation that maps $[8,3,4,5]$ to its minimal image $[1,2,3,6]$ under $G$. For example we can use the permutation
$(1,4,3,2,5,6,7,8)$. This maps the set of non-fixed points $\{1,2,6,7\}$ to the set
$\{4,5,7,8\}$. We then pick an ordering for the orbits in the minimal image, for example by length, and then break ties by smallest member of the orbit, and then map the list of orbits back. The ordering in the minimal image is therefore \([4,5,7,8]\).
Applying the inverse of our minimising element above, we obtain a full ordering of the vertices and hence an ordered partition: $[8|3|4|5|1|2|6|7]$. When we return from \textsc{Refine}, \textsc{Part} fixes every point, and hence we generate a minimal image under $G$ of these fixed points, with minimising permutation $(1,4,3,2,5,6,7,8)$. We used this permutation earlier -- this will not always happen, but since we use minimal image in many places, the same permutations often appear. The present node is a leaf, and we return.

Going back to the node where we split for the second time, there is another branch where we colour first 3, then 5, then 4. There we reach a leaf in the same way, with the minimising permutation $(1,7,4,5,3,2,8)$ mapping $[8,3,5,4]$ to $[1,2,3,5]$. The corresponding ordered partition is $[8|3|5|4|7|6|1|2]$.
Going back one step further, we see that the remaining two branches also split once more, and in total we have six leaves.
For completeness, we describe the partitions that belong to the remaining leaves:

If we go down the branch where we first fix  
 4, then 3, then 5, then the minimal image of $[8,4,3,5]$ under $G$ is $[1,2,3,5]$, with a minimising element $(1,7,8)(2,6,4)$. 
 The final partition into singleton cells, at the leaf, is
 $[8|4|3|5|6|2|1|7]$.
 Fixing first 4, then 5, then 3 gives minimal image $[ 1, 2, 3, 6 ]$ with minimising permutation $(1,4,2,7,5,3,6,8)$. The partition at the leaf is $[8|4|5|1|7|3|2|6]$.
For the branch where we fix first 5, then 3, then 4, we 
minimise $[8,5,3,4]$ with the element $(1,4,6,5,2,8)$, which gives the 
partition
$[8|5|3|4|1|6|7|2]$ at the leaf.
Finally, if we fix 5, 4 and 3 in this order, then 
the minimising element is $(1,7,6,8)(2,4,3,5)$ and the partition at the leaf is $[8|5|4|2|3|7|1|6]$. Therefore, in total, our search generates six minimising permutations and hence six candidates for a canonical image.
  \end{example}

Note that the execution of \Cref{alg-basic-search} is compatible with the action of $G$ on $\actSet$.
This is because each sub-function we execute during the algorithm is, by construction, compatible with the action of $\Sym\Omega$, except the refiner $R_G$. But this is compatible with the action of $G$.

\begin{remark}\label{rem-tree-compatible}
Suppose that \(g\in \Sym{\Omega}\), that \(a,b \in \actSet\) and that $a^g=b$. 
Moreover let $n,k \in \N$ and let \([R_1,...,R_n]\) be a list of compatible refiners for \(\actSet\), $[\Split_1,...,\Split_k]$ a list of compatible splitters and $\approxFunc$ a compatible approximator for producing a search tree for $a$ as in Algorithm \ref{alg-basic-search}.

Then mapping $T_a$ to \(T_a^g\) gives the same tree as producing $T_b$ by following Algorithm \ref{alg-basic-search}.
\end{remark}

We point out once more that, for this remark to hold, it is important that the children of nodes are treated as unordered in theory, although they will have an ordering in the implementation.

\section{Finding canonical images using search trees}\label{sec-images-from-trees}

In this section we describe how to find a canonical image from a search tree,
and we prove that what we find satisfies all conditions of a canonical image.
Rather than providing a mathematical description of the canonical image of an object \(a\), our proof will instead show directly that, if two objects \(a\) and \(b\) and are in the same orbit under some action, then they will have the same canonical image. Much of this proof will come from \textit{compatibility}, and analysing the relationship between the search trees for \(a\) and \(b\).

Another important requirement is that our algorithm can be implemented efficiently. We will mention briefly where the algorithm is designed to allow efficient implementation, at the cost of complicating the description. Further, at the end of this section we will discuss some improvements which greatly improve performance, without changing the canonical image.

\begin{remark}\label{lem-min-image-list}
Suppose that \(G \leq \Sym{\Omega}\) and that
$l$ is a list that contains each element of $\Omega$ exactly once.
Then there is a unique permutation, which we call a \textbf{minimising permutation}, in $G$ that maps $l$ to its minimal image under $G$ with respect to the lexicographic ordering on lists.
The minimal image of $l$ and the minimising permutation can be calculated in polynomial time (see appendix).
\end{remark}


\begin{example}\label{ex-lexi}
Let $\Omega\coloneqq\{1,...,8\}$ and 

$G\coloneqq
\langle(1,2,3,4,5,6,8), (1,3,2,6,4,5),(1,6)(2,3)(4,5)(7,8)\rangle \le \Sym\Omega$.\\
We consider the minimal images of the lists $[8,3,4,5,1,2,6,7]$ and $[8,3,5,4,1,2,6,7]$ with respect to lexicographical ordering of lists using the natural ordering of $\Omega$ -- under both $\Sym\Omega$ and $G$.
Since $\Sym\Omega$ is 8-transitive on $\Omega$, the minimal image of both lists
under $\Sym\Omega$ is \\
$[1,2,3,4,5,6,7,8]$, with minimising permutations $(1,5,4,3,2,6,7,8)$ and $(1,5,3,2,6,7,8)$ respectively.
Under $G$ we obtain the minimal images
$[1,2,3,6,4,5,7,8]$, with minimising permutation \\
$(1, 4, 3, 2, 5, 6, 7, 8)$, and 
$[ 1, 2, 3, 5, 7, 8, 6, 4 ]$ with minimising permutation $(1,7,4,5,3,2,8)$, respectively.
\end{example}

We will now explain the connection between our search trees and canonical image.

\begin{definition}\label{def-canonical}

Let \(G \leq \Sym\Omega\), \(\preceq\) denote a 
compatible partial order and let \partFunc{} denote an  ordered orbit approximator.

We define a function $\Can$ from \actSet to itself by following the steps below. We will prove later that this actually is a function, and that it maps members of \actSet to a canonical member of their orbit under $G$.  Let $a \in \actSet$ and let \(T_a\) be the search tree for $a$ as defined in \Cref{alg-basic-search}. Recall that this means that each node is labelled with a \emph{graph stack}, and leaf nodes are also labelled with a \emph{result} permutation.

  \begin{enumerate}
  \item   Denote the set of leaves of \(T_a\) as \(L_a\).
  \item For each leaf \(l \in L_a\) define a list of graph stacks \(\textsc{Stacks}(l)\) by taking, in order, the \emph{graph stack} vertex label of each vertex on the path from the root of \(T_a\) to \(l\). We remark that each member of this list will be a prefix of the next one, because graph stacks accumulate graphs as we go along a branch through the search tree.

  \item Extend the partial order \(\preceq\) to the elements of \(\{\textsc{Stacks}(l) \mid l \in L_a\}\) by
  ordering these lists lexicographically, comparing pairs of elements by \(\preceq\).
  
  \item Define $M_a \coloneqq \{m \in L_a \mid $ There is no $l \in L_a$ such that $\textsc{Stacks}(l) \prec \textsc{Stacks}(m)\}$. This means that $M_a$ consists of those elements of \(L_a\) which have minimal images with respect to \(\preceq\) when \textsc{Stacks} is applied.
 
  \item Build the set of permutations \(P_a = \{\textrm{result}(m)\ |\ m \in M_a\}\), where $\textrm{result}(m)$ is the result permutation label from the leaf node $m$.

  \item Finally, define \(\Can(a)\coloneqq\Min_{\actEq}\{a^g \mid g \in P_a\}\). 
  \end{enumerate}

\end{definition}

Note that the definition of \textsc{Stacks} in \Cref{def-canonical} as a list of graph stacks is chosen in this way to allow efficient algorithm implementation. The algorithm would still produce a correct canonical image if we only compared the graph stacks from the leaf nodes in steps 3 and 4 of \Cref{def-canonical}, and this would also simplify describing and proving the correctness of the algorithm. The advantage of using \textsc{Stacks} is that it allows us to sometimes prove a node is not in \(M_a\) in Step 4 by only considering a prefix  of \textsc{Stacks}.

\begin{lemma}\label{lem-canonical}

Let $\Can$ be as defined in Definition \ref{def-canonical}.
Then $\Can$ is a well-defined function from \actSet to itself, and moreover, it is a canonical labelling function with respect to the action of $G$ on \actSet.
\end{lemma}

\begin{proof}
Let $a \in \actSet$ and suppose that $T_a$ and $T'_a$ are search trees constructed as in Algorithm \ref{alg-basic-search}. Both trees start with the same root node, and the same refiners are applied in the same ordering.
A choice occurs only when a splitter is applied, because the children could appear in different orders for $T_a$ and $T'_a$.
This is relevant for implementation, but in theory the trees do not differ, so we obtain the same set $L_a$ of leaves in both cases.
This also means that the set 
\(\{\textsc{Nodes}(l) \mid l \in L_a\}\) 
is independent of whether we build the tree $T_a$ or the tree $T'_a$. Since we started with the same ordering \(\preceq\) and hence the same induced 
lexicographical ordering, it follows that 
the sets $M_a$ and $P_a$ are  also independent of the choice of the search tree.
This implies that \(\Can(a)\) is well-defined.

Next suppose that $b \in a^{G}$ and let $g \in G$ be such that $a^{g}=b$.
Then Remark \ref{rem-tree-compatible} tells us that
\(T_a^g\) is the same tree as if we produced $T_b$ following Algorithm \ref{alg-basic-search}.

Hence the set of leaves $L_b$ of $T_b$ is the same as $L_a^g$,
and for each $l \in L_b$ there is a leaf $t \in L_a$ such that
$t^g=l$, so $\textsc{Stacks}(l)=\textsc{Stacks}(t)^g$. Since $\preceq$ is compatible with applying $g$, we also know that $M_b = M_a^g$. Moreover, $\partFunc$ satisfies $\partFunc(a^g)=\partFunc(a)^g$ and \textsc{MinPerm} satisfies the condition that \(\textsc{MinPerm}(G,L^g) = \textsc{MinPerm}(G,L)^g\), and then it follows that \(P_a \cdot g=P_b\).

In our final step we consider the minimal members of \(P_a\) and \(P_b\):

\begin{align*}
 & \{a^{p^{-1}} | p \in P_a\} \\
=& \{a^{g{g^{-1}}{p^{-1}}} | p \in P_a\} \\
=& \{a^{g{{(pg)}^{-1}}} | p \in P_a\} \\
=& \{a^{g{q^{-1}}} | q g^{-1} \in P_a\} \textrm{ (where } q = pg )\\
=& \{a^{g{q^{-1}}} | q \in P_a \cdot g\}\\
=& \{b^{q^{-1}} | q \in P_b\}
\end{align*}

Therefore we have mapped \(a\) and \(b\) to a common set of images, namely \(\{a^{p^{-1}} | p \in P_a\}\) = \(\{b^{p^{-1}} | p \in P_b\}\).
Given the definition of the canonical image of \(x\) as \(min\{x^{p^{-1}} | p \in P_x\}\), for a fixed total ordering \(min\) on our objects, we deduce that $\Can(a)$ and $\Can(b)$ are equal.
\end{proof}

\begin{example}\label{ex-can-from-tree}
We continue the discussion from Example 
\ref{example:canonical-image-graph}.
There we had six leaves, and we consider Definition \ref{def-canonical} for each of them.
Starting with the leaf $l_1 \coloneqq [8,3,4,5,1,2,6,7]$, we follow the path from the root node to this leaf as follows:

We begin at the root node, with the empty graph stack. The first child has two graphs, created by the refiners, namely the graph $\Gamma$ itself from the refiner for $\Gamma$, and the complete graph $C_\Omega$ on $\Omega$ (from refining $G$ by orbits, because $G$ is transitive).

At this point we have completed the first node, and therefore, this will be the first member of \textsc{Stacks} for every leaf node. For the first child, we branch on $\Gamma_3$, $\Gamma_4$ and $\Gamma_5$ in turn, and in each branch we create a new node.  
Then we build a graph where vertices have the same colour if and only if they
are in the same orbit of $G_{[8,3,4,5]}$, which means that every vertex is individualised. This last graph will be abbreviated with $l_1$ itself, the singleton partition corresponding to it, because the ordering of the individualisation will matter in practice.
We conclude that $\textsc{Stacks}(l_1)$ looks like this:

\noindent
$[\varnothing,[\Gamma], 
[\Gamma, C_\Omega], 
[\Gamma, C_\Omega, \Gamma_3], 
[\Gamma, C_\Omega, \Gamma_3,\Gamma_4],[\Gamma, C_\Omega, \Gamma_3,\Gamma_4,\Gamma_5],$
\noindent$[\Gamma, C_\Omega,  \Gamma_3,\Gamma_4,\Gamma_5, l_1]]$,

The $i^{th}$ member of \textsc{Stacks} is the graph that was added to the stack at depth $i$, if we count from the root node. Note that in practice this is not stored as a list of stacks, but as a single stack, and the size of the stack at the end of each node. Also, we assume that refiners do not add duplicate graphs -- adding these would not change the results or the correctness of the algorithm (as long as the graphs were added consistently).

Using the same notation and looking at the remaining leaves in the ordering in which we describe them in Example \ref{example:canonical-image-graph}, the remaining lists coincide with 
$\textsc{Stacks}(l_1)$ up to the 3rd entry. Thus, we do not repeat these entries and only give the remaining ones.

\noindent$\textsc{Stacks}(l_2)=[
...,
[\Gamma, C_\Omega, \Gamma_3],
[\Gamma, C_\Omega,  \Gamma_3,\Gamma_5],[\Gamma, C_\Omega,  \Gamma_3,\Gamma_5,\Gamma_4],[\Gamma, C_\Omega,  \Gamma_3,\Gamma_5,\Gamma_4, l_2]]
$,

\noindent$\textsc{Stacks}(l_3)=[
...,
[\Gamma, C_\Omega, \Gamma_4],
[\Gamma, C_\Omega,  \Gamma_4,\Gamma_3],
[\Gamma, C_\Omega, \Gamma_4,\Gamma_3,\Gamma_5],
[\Gamma, C_\Omega, \Gamma_4,\Gamma_3,\Gamma_5, l_3]]
$,

\noindent$\textsc{Stacks}(l_4)=[
...,
[\Gamma, C_\Omega,  \Gamma_4],
[\Gamma, C_\Omega,  \Gamma_4,\Gamma_5],
[\Gamma, C_\Omega,  \Gamma_4,\Gamma_5,\Gamma_3],
[\Gamma, C_\Omega,\Gamma_4,\Gamma_5,\Gamma_3, l_4]]$,

\noindent$\textsc{Stacks}(l_5)=[
...,
[\Gamma, C_\Omega,  \Gamma_5],
[\Gamma, C_\Omega, \Gamma_5,\Gamma_3],
[\Gamma, C_\Omega,  \Gamma_5,\Gamma_3,\Gamma_4],
[\Gamma, C_\Omega,  \Gamma_5,\Gamma_3,\Gamma_4, l_5]]$,
and

\noindent
$\textsc{Stacks}(l_6)=[
...,
[\Gamma, C_\Omega,  \Gamma_5],
[\Gamma, C_\Omega, \Gamma_5,\Gamma_4],
[\Gamma, C_\Omega,\Gamma_5,\Gamma_4,\Gamma_3],$
\noindent $[\Gamma, C_\Omega, \Gamma_5,\Gamma_4,\Gamma_3, l_6]]$.

The final entries $l_2,...,l_6$ correspond to the remaining singleton partitions from Example 
\ref{example:canonical-image-graph}. 
In our ordering of graphs, the graphs 
$\Gamma_3$, \ldots, $\Gamma_8$ are considered to be the same because they have in common that only a single vertex is coloured and there are no arcs or edges.

We now attempt to order the different \textsc{Stacks} -- this is done using a permutation-invariant ordering, in lexicographic ordering. Looking at $\textsc{Stacks}(l_2)$ and $\textsc{Stacks}(l_3)$, the identity permutation maps the first three entries (which we omitted) of $\textsc{Stacks}(l_2)$ to the corresponding entries of $\textsc{Stacks}(l_3)$, the permutation $(3,4)$ maps the sixth entry to the sixth and $(3,4,5)$ maps the sixth and seventh entries simultaneously. This means that the search could not stop before reaching depth 8, because we cannot distinguish between the lists of graphs built at depths 1 to 7.

This means that in terms of the lexicographical ordering, the six stacks above only differ at the very last graph stack, and there only at the last entry.
These final graphs are written as lists, which means that they can be ordered lexicographically. We see that $l_1$ is the smallest one, therefore $M_a=\{l_1\}$ and $P_a=\{
(1,4,3,2,5,6,7,8)\}$.
Then the canonical image of the graph $\Gamma$ is defined to be its image under $(1,4,3,2,5,6,7,8)$, which is the graph with vertex set $\Omega$ and 
the following edges:
$\{4,7\}, \{4,8\}, \{5,7\}, \{5,8\}, \{2,1\}, \{3,1\}, \{6,1\}$.
\end{example}

\subsection{Improving Performance}

When running our algorithm in practice, we improve performance by skipping the generation of parts of the search tree when we can prove those parts of the search tree will not contain the canonical image. The techniques we use are similar to those used in other backtrack search algorithms, in particular Nauty and our previous work on Graph Backtracking.

The most significant technique involves skipping generating the children of nodes if we can prove that the leaves generated from this node are not elements of the set $M_a$ from Step 4 of \Cref{def-canonical}. Those parts of the tree which do not generate elements of $M_a$ do not need to be generated at all, during search, because they are not used. In practice, search trees are generated as a depth-first transversal, and we keep track of the current smallest list of graph stacks generated by any leaf. If some node of the tree produces a lexiocographically larger list of graph stacks than the smallest known list of graph stacks for a leaf, then all children of this node will also create lexicographically larger lists of graph stacks. Therefore, generation of the entire sub-tree under this node can be skipped, because the ordering of lists of graph stacks in \Cref{def-canonical} uses lexicographic ordering.

\section{Experiments}

In this section, we will measure the performance of our algorithm compared to alternatives. As there is no other implementation of a general framework for finding canonical images in arbitrary permutation groups, each experiment will, where appropriate, compare against specialised algorithms.
All of our experiments are performed with GAP 4.12~\cite{GAP4}. Our algorithm is implemented in
the Vole package \cite{vole}. For all instances, there is a time-out of 30 minutes and a memory limit of 4GB.
We note that we specifically search for difficult problems, both for our algorithm and
existing algorithms. 

\subsection{Grid group experiments}

As in our previous work on minimal and canonical images~\cite{minimal-canonical}, we will use ``grid groups''.
We use these because, in practice, they often produce difficult problems. We reproduce the
definition of a grid group here.

\begin{definition}\label{def:gridgroup}
  Let $n \in \N$. The direct product $\Sym{n} \times \Sym{n}$ acts on the set
  $\{1,\ldots, n\} \times \{1,\ldots, n\}$ of pairs in the following way:

  For all $(i,j) \in  \{1,\ldots, n\} \times \{1,\ldots, n\}$
  and all $(\sigma,\tau) \in \Sym{n} \times \Sym{n}$ we define

  \[
    {(i,j)}^{(\sigma,\tau)} \coloneqq (i^{\sigma}, j^{\tau}).
  \]

  The subgroup $G \le \Sym{\{1,\ldots, n\} \times \{1,\ldots, n\}}$ defined by this
  action is called the \textbf{$n \times n$ grid group}.
\end{definition}

We begin by reproducing the experiments of our minimal and canonical images paper \cite{minimal-canonical},
looking for the canonical image of randomly generated sets of points in an \(n \times n\) grid group.
Here, we consider sizes of grids from 5 to 50 inclusive (46 sizes of grids in total) and we note that a stabilizer chain for the group is calculated before starting the experiment. Otherwise the time is dominated by the time taking to build a stabilizer chain.

In \cite{minimal-canonical} we compared 16 search orderings. We consider the very best one (RareOrbitPlusRare),
and MinOrbit, as this most closely matches the search ordering Vole uses. Overall we considered
$138$ experiments, meaning that Vole solved every instance we considered within the allowed time limit. The
one experiment where RareOrbitPlusRare failed reached the memory limit of 4GB.
The table below gives the number of tests which were solved in the time and memory limit, and the total time (in seconds) that was taken by those experiments which finished. We note that this means that algorithms which only finished few tests may also have a low time.\\

\small
{
\begin{tabular}{|l|r|r|r|r|r|r|r|r|}
\hline
Size of set&\multicolumn{2}{|c|}{\(\myfloor{\frac{n^2}{2}}\)}&\multicolumn{2}{|c|}{\(\myfloor{\frac{n^2}{4}}\)}&\multicolumn{2}{|c|}{\(\myfloor{\frac{n^2}{8}}\)}&\multicolumn{2}{|c|}{total} \\
Results &solved &time&solved &time&solved &time&solved &time\\
\hline
Vole & 46 & 53& 46 & 58& 46 & 106& 138 & 218\\
RareOrbitPlusRare & 45 & 185& 46 & 18& 46 & 10& 137 & 214\\
MinOrbit & 8 & 60& 13 & 260& 21 & 116& 42 & 438\\
\hline
\end{tabular}}\\

\normalsize
\smallskip
Vole is overall the same speed as the RareOrbitPlusRare algorithm. This is
despite the fact that it contains no specialised code for finding canonical images of sets of points and can also be used to find canonical images of more complex combinatorial structures, such as sets of sets, or sets of tuples, which are not supported by the algorithm in \cite{minimal-canonical}.

Next we looked at canonical images of permutations under conjugation in the grid group. Since there is no other tool available for this purpose, as far as we know, we compared
the time it took to find the canonical image of two permutations, $g$ and $h$, under a group $G$ against the time it takes GAP to find a
permutation which can map $g$ to $h$ in $G$, known as \texttt{RepresentativeAction}. Since the canonical images of $g$ and $h$
can be used to implement \texttt{RepresentativeAction}, this is a reasonable comparison.

However, this can be viewed as slightly unfair, with an advantage for Vole. If we wanted to sort $n$ permutations into conjugacy classes, for example, then we would only have to call Vole $n$ times, whereas GAP may require calling \texttt{RepresentationAction} on all pairs of permutations. 

We chose to only use permutations built from pair-wise disjoint 2-cycles (i.e. involutions). When the size of cycles grows, and when there are cycles of many different sizes, then both Vole and GAP solve problems very quickly and the time taken becomes dominated by the time taken to build an initial stabiliser chain for the group. For each size of grid, we consider ten randomly generated permutations. The results show both the number of problems solved and the average time taken for the solved instances in seconds.\\

\begin{tabular}{|l|r|r|r|r|r|r|r|r}
\hline
$n$&\multicolumn{2}{|c|}{GAP}&\multicolumn{2}{|c|}{Vole}\\
&\# solved&time (sec)&\# solved&time (sec)\\
\hline

10 &10 & 0.9 &10 & 0.7 \\
20 &1 & 58 &10 & 13 \\
30 &0 & - &10 & 101 \\
40 &0 & - &10 & 448 \\
\hline
\end{tabular}

\vspace{0.5cm}
We can see that Vole is much faster than GAP and able to solve problems far beyond the ability of GAP.

Out of curiosity, we also looked at the performance of \textsc{Magma} (see \cite{Magma}) in checking involutions in grid groups for conjugacy. 
\textsc{Magma} solved this problem much faster than either GAP, or Vole. It solved the problems in $50 \times 50$ grid groups in 0.5 to 7 seconds. According to the documentation for \textsc{Magma}, this functionality is implemented using Leon's Partition Backtrack~\cite{leon1997}, while greatly outperforming the implementations in GAP and Vole -- at least when we only need to check if two permutations are conjugate. 
This is a constant reminder
for us that further improvements are possible and necessary.
We were unable to find out why \textsc{Magma} is able to solve these problems so quickly, compared to the implementation of Leon's Partition backtrack code in GAP. For example, it would be interesting to know whether the impressive speed is possible because of the \textsc{Magma} function 
\textsc{ClassRepresentative}. 
If that was used and was at least  partly responsible for the performance, then our work on canonical images might pave the way for future speed-ups in \textsc{GAP}, and in other systems as well. 

\subsection{Canonical images under the action of primitive groups}

In this experiment we reproduce an experiment from \cite{minimal-canonical}, where we look for canonical images of randomly generated sets under the action
of a primitive group which moves between \(10\) and \(300\) points. We remove the natural alternating and symmetric groups, because these groups can easily be treated as special cases.

The results are given below, showing both the best results for static and dynamic selection strategies from \cite{minimal-canonical}, and Vole. The total number of instances solved, and the time in seconds taken for these solvable problems is given. A total of 5,952 experiments were attempted, meaning no system completed every problem. 

\tiny{
\begin{tabular}{|l|r|r|r|r|r|r|r|r|}
\hline
Size of sets
&\multicolumn{2}{|c|}{\(\myfloor{\frac{n}{2}}\)}&\multicolumn{2}{|c|}{\(\myfloor{\frac{n}{4}}\)}&\multicolumn{2}{|c|}{\(\myfloor{\frac{n}{8}}\)}&\multicolumn{2}{|c|}{total} \\
Search&\# solved&time (sec)&\# solved&time (sec)&\# solved&time (sec)&\# solved&time (sec)\\
\hline
Vole & 1963 & 7517& 1984 & 8573& 1984 & 2410& 5931 & 18501\\
RareOrbitPlusRare & 1956 & 2653& 1983 & 627& 1984 & 46& 5923 & 3327\\
FixedMinOrbit & 1916 & 4139& 1979 & 907& 1984 & 69& 5879 & 5116\\
\hline
\end{tabular}}\\

\normalsize

While Vole does take longer on average, it solves a larger number of instances. The most significant reason for the slowdown is that Vole is generating orbital graphs for the groups as search progresses. In some cases these graphs are very large and do not reduce the size of the search. In the future we plan to look for better heuristics to filter which graphs are generated and used during refinement. Also, in this experiment, we did not re-use the graphs on different calls to Vole.

\subsection{Canonical images of graphs}

Vole can solve a much wider range of problems than previous graph canonising systems such as Nauty/Traces \cite{nauty} and Bliss \cite{bliss}. This is because it supports finding canonical images under any group of symmetries on the vertices of the graph, while Nauty, Traces and Bliss can only search within the full symmetric group on the vertex set.

Here, we therefore compare Vole and Nauty on the canonical image problems that Nauty can solve. We use the latest version of Nauty, and we found that this was faster than Bliss for these experiments.

On these problems, there is no fundamental reason that Vole should not perform similarly to Nauty. However, in practice, it is significantly slower in terms of the fastest speed Vole can be, a slowdown of around 5 times.
After investigating how our algorithm performed, we believe that there are three main reasons for this:

\begin{itemize}
    \item  Vole supports refining stacks of graphs, and the code is less efficient when only a single graph is refined. 
    \item Vole uses
stabiliser chains to build the known stabiliser as the search progresses -- this calculation is done in GAP and therefore requires communication costs. 
\item Nauty uses a different, better ordering to decide which cell to split on.
\end{itemize}

In general, Nauty has been carefully optimised for over 40 years, while Vole is comparatively new, and hence there are many
small optimisations which still have to be considered and implemented. We believe there is no optimisation in
Nauty which could not, in principle, be added to Vole.

Vole can solve graph problems which Nauty cannot -- finding canonical images under the action of groups other than the symmetric
group on the set of vertices. For example, we look for the canonical image of graphs under the action of cyclic groups and alternating groups on the whole set of vertices. In our experiment, we compare Nauty and Vole in the symmetric group, because Nauty can only compute canonical images of graphs under the full symmetric group. We also use Vole to search for canonical images of graphs in the cyclic group and alternating group on the set of vertices -- showing how Vole can solve interesting problems that have been impossible to solve with previous graph automorphism systems such as Nauty.

The slowdown that we observe when we search for canonical images under subgroups, rather than under the full symmetric group, is entirely due to calculating stabiliser chains. The only reason why it does not occur for the symmetric groups is that they are treated as a special case.  In our experiments, the searches in smaller groups explore the same number of nodes, or even fewer nodes. This behaviour is not guaranteed, but it usually occurs in practice. We could add similar optimisations for the cyclic and alternating groups, but in general our algorithm's performance is greatly affected by the performance of stabiliser chain calculations, which is why we plan to improve these calculations in the future. In particular, currently GAP, and therefore Vole which re-uses GAP stabiliser chain implementation, cache at most a single stabiliser chain at a time, producing repeated work if a previously created stabiliser chain is needed again.

For our experiments we use the ``Steiner Triple System Graphs'' from the website of Nauty and Traces, which were originally distributed with Bliss 
(\url{https://pallini.di.uniroma1.it/Graphs.html}). The times in these experiments are given in seconds.\\

\begin{tabular}{|l|r|r|r|r|}
\hline
instance & nauty(symmetric) & Vole (symmetric) & Vole (cyclic) & Vole (alternating) \\
\hline
sts-13 & 0 & 0.014 & 0.055& 0.099 \\
sts-15 & 0 & 0.016 & 0.071& 0.309 \\
sts-19 & 0.001 & 0.031 & 0.134& 1.235 \\
sts-21 & 0 & 0.034 & 0.184& 1.994 \\
sts-25 & 0.009 & 0.086 & 0.329& 9.398 \\
sts-27 & 0.001 & 0.099 & 0.423& 18.848 \\
sts-31 & 0.035 & 0.251 & 0.869& 51.666 \\
sts-33 & 0.001 & 0.233 & 1.019& 133.580 \\
sts-37 & 0.089 & 0.601 & 1.441& 265.022 \\
sts-39 & 0.004 & 0.632 & 1.953& 623.571 \\
sts-43 & 0.219 & 1.409 & 2.918& 1124.527 \\
\hline
\end{tabular}

\subsection{Results overview}

Overall, we can see that Vole has the potential to become competitive with previous work on finding canonical images, after some optimization, with the exception of Nauty. Also, Vole can solve a wide range of other canonical image problems which were previously insolvable by any system, and these can be solved very efficiently.

\section{Conclusion}

In this paper we present a general framework for finding canonical images, and demonstrate its practical value on a range of experiments. Compared to previous systems, our implementation Vole can solve a wider range of different canonical image problems.

One significant advantage of our framework is that the same refiners can often be used to solve intersection and coset problems. We hope this will create a virtuous circle, as new refiners will have, in many cases, immediate applications to both intersection and canonical image problems. Almost all refiners in this paper were re-used, unchanged, from our earlier work -- only one new refiner (described in Appendix 7.1) was created, and that refiner is now the refiner Vole uses for a group given by a list of generators. In the future we plan to work on refiners for normalisers, which based on the framework in this paper will also be used for finding canonical images of groups under conjugacy.

\section{Appendix}

Here we present two algorithms which are needed in this paper. Firstly, we present how to find the minimal image of a list of elements of \(\Omega\) with respect to the action of a group \(G \leq \Sym{\Omega}\). Secondly, we present a new refiner for a group given as a list of generators. Previously presented refiners for a group given as a list of generators could return the cells of partitions in different orders in different searches. This was not an issue previously, but it is for the algorithm we present in this paper.

In this section we say that a list of elements of $\Omega$ is \textbf{exhaustive and non-repetitive} if it contains each element of $\Omega$ exactly once.

\begin{lemma}
Given a set \(\Omega\) with a total ordering \(\leq\), and a group \(G \leq \Sym{\Omega}\), let \(\textsc{MinPerm}_G\) denote the function 
that maps each non-repetitive exhaustive list on $\Omega$ to its unique minimising permutation in $G$.

Then the images under \(\textsc{MinPerm}_G\)
 can be calculated in polynomial time. 
\end{lemma}

\begin{proof}
If we have two exhaustive, non-repetitive lists on $\Omega$, then there is a unique permutation in $\Sym\Omega$ that maps one to the other.
Therefore \(\textsc{MinPerm}_G\) is a well-defined function. 

Given a list \(L\) with members drawn from \(\Omega\), the minimal image of \(L\) under \(G\) can be found in a series of steps. Firstly, the minimal image under a lexicographic ordering  must map the first member \(L[1]\) of \(L\) to its smallest value \(m_1\) under \(G\) with respect to the ordering $\leq$ on $\Omega$. Finding a permutation \(g \in G\) which maps \(L[1]\) to \(m\) can be done by finding the orbit of \(L[1]\) under \(G\). The minimal image of \(L\) under \(G\) is the same as the minimal image of \(L^g\) under the point stabilizer \(G_{m_1}\), because in the minimal image of \(L^g\), \(m_1\) must be fixed. We can continue mapping each element of \(L\) to its minimal image under the subgroup of \(G\) which maps all previous members of \(L\) to their smallest values, respectively. This calculation is finished after finitely many steps. Calculating orbits and point stabilizers in permutation groups can be performed in polynomial time using stabiliser chains, and consequently the whole calculation can be done in polynomial time.
\end{proof}

\subsection{Refiners for groups defined by minimal images}

Searching for canonical images inside a group $G$ requires a refiner for $G$. Typically $G$ will be specified as a subgroup of a larger permutation group, with a set of generators. 
A family of refiners for a group given by a set of generators is given in Section 7.3.2 of \cite{gb-extended}, and it uses the orbital graphs of the group.
This
extends a previous refiner given in \cite{newrefiners}. Leon also gives a refiner which uses only the orbits of the group \cite{leon1997}. These algorithms are different to any other refiners presented in the literature, as they do not describe a single refiner -- instead they describe how to build a refiner dynamically during search. This means that the refiners can produce different outputs when solving different problems.

These refiners are not suitable for use in our canonical image algorithm, because the proof of correctness requires that a single fixed refiner is used for $G$. In this section we present a new refiner for groups given by a list of generators. The output of this refiner is almost identical to the previously presented refiner, the only difference is that the output is fixed, and therefore it produces the same output in every search, given the same input graph stack.

We begin by discussing the most important pieces of our refiner, and then we use this to describe the problem previous refiners were trying to solve and to explain our new solution for this problem.

\begin{definition}\label{def-fixed-refiner}
Given a group $G$ on a set $\Omega$, a \textbf{fixed point refiner} $F$ for $G$ is a function from the set of lists of points of $\Omega$ to $\Stacks{\Omega}$ such that for each list $L$ of points from $\Omega$ and \(g \in G\)
the following holds: If \(L^g=L\), then \(F(L)^g=F(L)\).
\end{definition}

Note that a fixed point refiner, as given in \Cref{def-fixed-refiner}, does not define a valid refiner in the sense of 
Definition \ref{defn-comp-refiner} and \cite{gb-published} for two reasons -- its input is a list of fixed points, rather than a graph stack, and it defines only a single function, rather than a pair of functions.

\begin{example}\label{ex-fixed-refiner}
Given a group $G$ on a totally ordered set $\Omega$, we give two examples for fixed point refiners. In both cases we make use of $G_L$, the point-wise stabilizer in $G$ of a list $L$ of points of $\Omega$.

The Orbit refiner \textsc{FOrbit} maps a list $L$ to a graph on $\Omega$ with no edges, where the vertices are labelled with the smallest member of $\Omega$ in their orbit under $G_L$.

The Orbital refiner \textsc{FOrbital} maps a list $L$ to a graph stack, where each graph is an orbital graph of $G_L$. These graphs are ordered in lexicographical order, using the total ordering of $\Omega$.
\end{example}

The orbit refiner from \Cref{ex-fixed-refiner} is similar to the technique used by Leon \cite{leon1997}, while the orbital refiner is the technique used in \cite{newrefiners} and \cite{gb-extended}. The most obvious way to turn a fixed point refiner into a refiner in the sense of \cite{gb-published} would be to take the fixed points of the graph stack using a fixed point approximator, and then apply the fixed point refiner to its output. However, this obvious way does not work: Given a refiner $(R,R)$ for a group $G$ (refiners for groups always consist of a pair of identical functions), it must satisfy the condition that for all $g \in G$ and all graph stacks $S$, it is true that $R(S^g)=R(S)^g$. The way we order the orbits and orbital graphs in \Cref{ex-fixed-refiner} will not satisfy this, because the orbits and orbital graphs are sorted. In previous work this problem was solved by storing all the orbits, or orbital graphs, ever created. Then, if a list of fixed points $L'$ was created which could be mapped to a previously created list $L$ by some $g \in G$, then the image of the refiner output for $L$ under $g^{-1}$ could be returned.

The reason for the dynamic nature of this refiner is that it takes the first list of fixed points from each orbit as a ``canonical member'' of that orbit. We deal with this problem by explicitly choosing the ``canonical member'' we will use for each orbit in advance -- the minimal image of the list of fixed points.

\begin{lemma}\label{lem-new-refiner}
Given a fixed point refiner $F$ for a subgroup $G$ of $\Sym\Omega$ and a fixed point approximator $\fixedFunc$, we define the function $\Refi$ from $\Stacks{\Omega}$ to itself as follows:

If $S \in \Stacks{\Omega}$, $L\coloneqq\Fixed{S}$ and $g \in G$ is a permutation which maps $L$ to its minimal image under $G$, then $\Refi(S)\coloneqq F(L^g)^{g^{-1}}$.

Then $(\Refi,\Refi)$ is a refiner for $G$ in the sense of \cite{gb-published}.
\end{lemma}

\begin{proof}

Firstly, we prove that $\Refi$ is well-defined, as the definition defines $g$ as any permutation in $G$ which maps $L$ to its minimal image under $G$. 
For this let $L$ be a list of points in $\Omega$ and let $g,h \in G$ be such that $M\coloneqq L^g=L^h$ is the minimal image of $L$ under $G$. Then $M^{g^{-1}h} = M$, so by the definition of a fixed point approximator, $F(M)^{g^{-1}h}=F(M)$. This means that $F(M)^{g^{-1}} = F(M)^{h^{-1}}$ and therefore $F(L^g)^{g^{-1}} = F(L^h)^{h^{-1}}$.

Next we show that $\Refi$ is a valid refiner for $G$ as explained in Definition 4.1 of \cite{gb-published}.  For all $\stackS, \stackT \in \Stacks{\Omega}$ such that
  $\stackS \cong \stackT$, we must prove that
  \[
    G \cap \Iso{\stackS}{\stackT}
    \subseteq
    G \cap \Iso{\Refi(\stackS)}{\Refi(\stackT)}.
  \]
  
It suffices to show that for all \(g \in G\), where $g \in \Iso{\stackS}{\stackT}$, $g$ is also in $\Iso{\Refi(\stackS)}{\Refi(\stackT)}$.

Let $L_S \coloneqq \Fixed{\stackS}$ and $L_T \coloneqq \Fixed{\stackT}$. By the definition of a fixed point approximator, the fact that $\stackS^g=\stackT$ implies that $L_S^g=L_T$. Therefore $L_S$ and $L_T$ are in the same $G$-orbit, and therefore they have the same minimal image $L_M$. 
Further, let $h \in G$ be such that $L_T^h=L_M$ and note that this means that $gh$ maps $L_S$ to $L_M$.

Then $\Refi(\stackS) = F(L_S^{gh})^{(gh)^{-1}} = F(L_M)^{(gh)^{-1}} = (F(L_M)^{h^{-1}})^{g^{-1}}$ and $\Refi(\stackT) = F(L_T^h)^{h^{-1}} = F(L_M)^{h^{-1}}$. This shows that $g \in \Iso{\Refi(\stackS)}{\Refi(\stackT)}$, completing the proof.
\end{proof}

In practice, implementing the refiner in \Cref{lem-new-refiner} is very similar to implementing the refiners for a group given by a list of generators as previously described in \cite{newrefiners, gb-extended,leon1997}. They each involve finding permutations which map a list of fixed points to a ``standard list'' -- the difference is that, in the previous algorithms, this standard list was chosen as the first element of each orbit to occur, while now it is fixed as the minimal image of each orbit. The old and new refiners will, when using the same general technique (either orbits, or orbital graphs), produce the same graph stacks, up to re-ordering of the graphs and re-naming the labels used for vertices and edges. We benchmarked this new algorithm against implementations of the previous refiners and found no measurable difference.

\subsection*{Acknowledgments}\label{sec-acknowledgments}

Many of the ideas that lead to this work 
stem from our project on Graph Backtracking, which is why we thank 
the VolkswagenStiftung (\textbf{Grant no.~93764})
and the Royal Society (\textbf{Grant code URF\textbackslash R\textbackslash 180015}) for their financial support of various aspects of this project.
The more recent work has been supported by
 the DFG (\textbf{Grant no.~WA 3089/9-1}) and again the Royal Society
(\textbf{Grant codes} \textbf{RGF\textbackslash EA\textbackslash 181005} and \textbf{URF\textbackslash R\textbackslash 180015}). We are grateful to all these institutions for their support.
Also, we would like to thank the referees for many valuable comments and suggestions on an earlier version of this article.

  \end{document}